\documentclass[11pt]{article}
\usepackage{amsmath,latexsym,amsfonts,amssymb,mathrsfs}
\usepackage{amsbsy}
\usepackage[normalem]{ulem} 
\newcommand{\stkout}[1]{\ifmmode\text{\sout{\ensuremafth{#1}}}\else\sout{#1}\fi}

\usepackage[usenames]{color}

\newcommand{\W}{{\mathcal W}}

\usepackage{tikz}

\newcommand{\J}[2]{\mathcal{J}_{#1}^{#2}}
\newcommand{\XT}[1]{X_{T,#1}}

\newcommand{\I}{\mathcal{I}}

\newcommand{\bsx}{{\boldsymbol{x}}}

\numberwithin{equation}{section}

\newcommand{\iprod}[1]{\langle#1 \rangle}
\newcommand{\biggiprod}[1]{\Big\langle#1 \Big\rangle}
\newtheorem{theorem}{Theorem}[section]

\newtheorem{lemma}{Lemma}[section]
\newtheorem{corollary}{Corollary}[section]

\title{An $\alpha$-robust and  second-order accurate scheme for a 
subdiffusion equation\thanks{School of Mathematics and Statistics, University of New South Wales, Sydney, Australia \\ This work was supported by the Australian Research Council grant DP220101811.}}
\author{Kassem Mustapha,  William McLean and Josef Dick}
\begin{document}

\maketitle
\begin{abstract}
We investigate a second-order accurate time-stepping scheme for solving  a time-fractional diffusion equation with a Caputo derivative  of  order~$\alpha \in (0,1)$. The basic idea of our scheme is based on local  integration followed by linear interpolation. It reduces to the standard Crank--Nicolson scheme in the classical diffusion case, that is, as $\alpha\to 1$. Using a novel approach, we show that the proposed scheme is $\alpha$-robust and second-order accurate in the $L^2(L^2)$-norm, assuming a suitable time-graded  mesh. For completeness, we use the Galerkin finite element method for the  spatial  discretization and discuss the error analysis under reasonable regularity assumptions on the given data.  Some numerical results are presented at the end.  
\end{abstract}
\section{Introduction}
We shall approximate the solution of the time-fractional diffusion equation 
\begin{equation}\label{eq:L1}
\partial_t^{\alpha} u(\bsx,t)+\mathcal{A} u(\bsx,t)=f(\bsx,t)
\quad\text{for $(\bsx,t)\in \Omega\times (0,T]$,}
\end{equation}
subject to homogeneous Dirichlet boundary conditions, that is, $u(\bsx,t)= 0$ 
on~$\partial\Omega\times (0,T]$, with $u(\bsx,0)=u_0(\bsx)$ at the initial 
time level~$t=0.$ The spatial domain~$\Omega\subset \mathbb{R}^d$ 
(with $d=1$, $2$, $3$) is a convex polyhedron, $0<\alpha<1$, the time 
fractional Caputo derivative 
\[\partial_t^{\alpha} v(t):=\I^{1-\alpha}v'(t)
=\int_0^t\omega_{1-\alpha}(t-s)v'(s)\,ds,\quad{\rm  with}\quad \omega_{1-\alpha}(t):=\frac{t^{-\alpha}}{\Gamma(1-\alpha)},\]
where  $v'=\partial v/\partial t$~and $\Gamma$ denotes the gamma  function. We use the notation $\I v(t)$ for the standard time integral of~$v$ from~$0$ to~$t.$ In \eqref{eq:L1}, ${\mathcal A}$ is  an elliptic operator in the spatial variables, defined  by~${\mathcal A}w(\bsx)=-\nabla\cdot(\kappa(\bsx)\nabla w)(\bsx)$. The  diffusivity~$\kappa \in L^\infty(\Omega)$ satisfies 
$0< \kappa_{\min} \le \kappa$ on $\Omega,$  for some constant~$\kappa_{\min}$. For the error analysis, we also require that 
$\kappa\in W^{1,\infty}(\Omega)$.

The presence of the nonlocal time fractional (Caputo) derivative 
in~\eqref{eq:L1} and the fact that the solution~$u$ suffers from a weak 
singularity near~$t=0$ have a direct impact on the accuracy, and consequently 
the convergence rates, of numerical methods.  To overcome this difficulty,
different approaches have been applied including corrections, graded meshes, 
and convolution quadrature \cite{GunzburgerWang2019,JinZhou2017,Kopteva2021,Mustapha2020,StynesOriordanGracia2015,WangWangYin2021}. Indeed, the numerical 
solutions for model problems of the form~\eqref{eq:L1}, including 
\emph{a priori} and \emph{a posteriori} error analyses and fast algorithms, 
were studied by various authors over the past fifteen years using multiple 
approaches \cite{Alikhanov2015,BanjaiMakridakis2022,CenHuangLeXu2018,HouHasanXu2018,JinLiZhou2018,JinLiZhou2020,Karaa2018,Kopteva2020,MaskariKaraa2019,ShenSheng2019,WangYanYang2020,WuZhou2021,YanEgwuLiang2021,ZengLiLiuTurner2015,ZhuXu2019}. For more references and details, see the 
recent monograph by Jin and Zhou~\cite{JinZhou2023}. 

In this work, we investigate rigorously the error from approximating  the solution of the initial-boundary value problem~\eqref{eq:L1} using a uniform second-order accurate time-stepping method. The latter is defined via  a local time-integration of problem \eqref{eq:L1} on each subinterval of the time mesh combined with continuous piecewise linear interpolation. The proposed scheme is identical to the piecewise-linear case of a discontinuous  Petrov--Galerkin method proposed in ~\cite{MustaphaAbdallahFurati2014}. Therein, with $\tau$ being the maximum time mesh step size,  a suboptimal convergence rate of order $O(\tau^{(3-\alpha)/2})$ was proved.  A time-graded  mesh~\eqref{eq: time mesh} was employed to compensate for the singular behaviour of the continuous solution at~$t=0$. In the limiting case as~$\alpha \to 1$, the problem~\eqref{eq:L1} reduces to the classical diffusion equation, and the considered numerical scheme reduces to the classical Crank--Nicolson method. In this case, $O(\tau^{(3-\alpha)/2})$ reduces to $O(\tau)$ which  is far from the optimal $O(\tau^2)$ rate achieved in practice. 

By using an innovative  approach that relies  on interesting implicit polynomial interpolations and  duality arguments, we show   $O(\tau^2)$ convergence, whilst at the same time relaxing the imposed regularity assumptions from the earlier 
analysis~\cite{MustaphaAbdallahFurati2014}. This convergence rate is $\alpha$-robust in the sense that the constant in the error bound remains bounded as $\alpha\to 1$. Implementation wise, although the proposed scheme is uniformly second-order accurate, the computational cost is comparable to the well-known backward Euler or $L1$~\cite{LiaoLiZhang2018,Mustapha2020,Stynes2022} methods, which are not even first-order accurate.  

For completeness, we discretize the problem~\eqref{eq:L1} over the spatial domain~$\Omega$ using  the standard Galerkin finite element method (FEM), thereby defining a fully discrete approximation to~$u$. An additional error of order $O(h^2)$ is anticipated under certain regularity assumptions on the continuous solution, where $h$ is the maximum spatial finite element mesh size. This is proved  via a concise approach that relies on the discrete version of the earlier error analysis. To make this feasible, the solution of the semidiscrete Galerkin finite element solution of  problem~\eqref{eq:L1} plays the role of the comparison function.

Outline of the paper. In the next section, we define our 
time-stepping scheme, introduce some notations and technical lemmas, and 
summarize the convergence results in Theorem~\ref{time-convergence}.
The required regularity properties are also highlighted. Section~\ref{sec:average}
proves some error bounds for the implicit piecewise-linear interpolant~$\widehat u$ 
defined in~\eqref{widehatu}.  Section~\ref{sec:semi} is devoted to showing 
the second order of accuracy of the proposed time-stepping scheme via a duality argument. 
In Section~\ref{sec:fully}, we discretize in space via the Galerkin finite
element method and discuss the convergence of the fully discrete solution.
To support our theoretical findings, we present some numerical results in 
Section~\ref{Sec: Numeric}.  Finally, a short technical appendix derives an $\alpha$-robust interpolation estimate.

\section{Time-stepping scheme}\label{Sec: FEM}
This section is devoted to discretizing the model problem~\eqref{eq:L1} over 
the time interval $[0,T]$ through a second-order accurate method, and stating
our main convergence results.  We begin by introducing some notations. 

For $\ell\ge 0,$ the norm on~$H^\ell(\Omega)$ is denoted by~$\|\cdot\|_\ell$.   
The Sobolev spaces $H^\ell(\Omega)$ and $H^1_0(\Omega)$ are defined as usual,
and the norm~$\|\cdot\|_{\dot H^r(\Omega)}$ in the (fractional-order) Sobolev 
space $\dot H^r(\Omega)$ is defined in the usual way via the Dirichlet 
eigenfunctions of the self-adjoint elliptic operator~$\mathcal{A}$ on~$\Omega$.
The inner product in~$L^2(\Omega)$ is denoted by $\iprod{\cdot,\cdot}$, and the
associated norm by~$\|\cdot\|$. The generic constant~$C$ remains bounded for $0<\alpha \le 1,$ and is independent of the time mesh and the finite element mesh, but may depend on
$\Omega$, $T$, and other quantities, including $\kappa$, $u_0$~and $f$.
 
Define the time mesh $0=t_0<t_1<t_2<\ldots<t_N=T$ by
\begin{equation} \label{eq: time mesh} 
t_n=(n\,\tau)^\gamma,\quad\text{with $\tau=T^{1/\gamma}/N$ and $\gamma \ge 1$,}
\quad\text{for $0\le n\le N$,}
\end{equation}
and let $\tau_n = t_n-t_{n-1}$.  Such a time-graded mesh has the
properties~\cite{McLeanMustapha2007}
\begin{equation}\label{time mesh properties}
t_n\le 2^\gamma t_{n-1} \quad\text{and}\quad 
\gamma \tau t_{n-1}^{1-1/\gamma}\le \tau_n\le \gamma \tau t_n^{1-1/\gamma},
\quad\text{for $n\ge 2$.}
\end{equation} 
Integrating problem~\eqref{eq:L1} over~$I_n:=(t_{n-1},t_n)$ and then dividing by~$\tau_n$
yields
\begin{equation}\label{eq: ibvp int}
\frac{1}{\tau_n}\int_{I_n} \partial_t^\alpha u\,dt 
	+\mathcal{A}\bar u_n=\bar f_n,\quad {\rm for}~~1\le n\le N,
\end{equation}
where $\bar f_n=\tau_n^{-1}\int_{I_n}f(t)\,dt$ denotes the average value of a function~$f$ over the time interval~$I_n$. Motivated by~\eqref{eq: ibvp int}, for $t\in I_n$  and for $1\le n\le N,$ our semidiscrete approximate solution $U(t)\approx u(t$  is defined by requiring that 
\[U(t)=\frac{t_n-t}{\tau_n} U^{n-1}
	+\frac{t-t_{n-1}}{\tau_n} U^n, \quad U^n:=U(t_n),\]
with \begin{equation} \label{fully}
\frac{1}{\tau_n}\int_{I_n}\partial_t^\alpha U\,dt
+\mathcal{A} U^{n-1/2}=\bar f_n,\quad{\rm with}~~U^0=U(0)=u_0,
\end{equation}
where $U^{n-1/2}=\bar U_n =\tfrac12( U^n+ U^{n-1}).$ 
At the $n$th time step, with $\tau_{n,\alpha}=\Gamma(3-\alpha)\tau_n^\alpha$, we have to solve an elliptic problem
for $W^n=U^n-U^{n-1}$, namely,
\[
(I+\tfrac12\tau_{n,\alpha}\mathcal{A})W^n=\tau_{n,\alpha}\bar f_n
    -\tau_{n,\alpha}\mathcal{A} U^{n-1}-\tau_{n,\alpha}\sum_{j=1}^{n-1}(\tau_n\tau_j)^{-1}b_{n,j}W^j,
\]
where  $b_{n,j}=a_{n,j}-a_{n,j+1},$ with 
$a_{n,j} =\omega_{3-\alpha}(t_n-t_{j-1})-\omega_{3-\alpha}(t_{n-1}-t_{j-1}).$


If $\alpha\to1$, then
$\partial_t^\alpha u\to u'$ and $\partial_t^\alpha U\to U'$, implying that our scheme reduces to the Crank--Nicolson scheme for the classical diffusion equation.

Our convergence analysis  relies on decomposing the error
as
\begin{equation}\label{eq:eta-decomposition}
\eta=u- U=\psi-\theta\quad\text{with}\quad\psi=u-\widehat u
\quad\text{and}\quad\theta= U-\widehat u,
\end{equation}
where $\widehat u$ is a continuous piecewise-linear function in time satisfying
\begin{equation}\label{widehatu}
\int_{I_n}\widehat u(t) \,dt=\int_{I_n} u(t)\,dt
\quad\text{for $1\le n\le N$,}\quad
\text{with $\widehat u(0)=\widehat u^0=u_0$.}
\end{equation}
Alternatively, $\widehat u$ can be defined via $\I\widehat u(t_n)=\I u(t_n)$ for $1\le n\le N$, with $\widehat u(0)=u_0$, and we say that
$\widehat u$ interpolates $u$ \emph{implicitly}.   The
decomposition~\eqref{eq:eta-decomposition} of the error~$\eta$ follows a
well-known pattern, but the novel choice of the piecewise linear
function~$\widehat u$ makes possible our improved error analysis under reasonable regularity assumption. 
The continuous average of~$u$ equals both the continuous and the discrete
average of~$\widehat u$ on each time subinterval~$I_n$.  For comparison,
let $u_I$ denote the usual continuous piecewise-linear interpolant to~$u$, that
is,
\begin{equation}\label{eq:uI}
u_I(t)=\frac{t_n-t}{\tau_n}\,u(t_{n-1})+\frac{t-t_{n-1}}{\tau_n}u(t_n)
\quad\text{for $t\in I_n$,}
\end{equation}
and observe that $u_I$~and $u$ have the same \emph{discrete}
average~$\tfrac12(u(t_n)+u(t_{n-1}))$ on each~$I_n$, but their continuous averages
will differ unless $u$ is linear on~$I_n$. 

Subtracting \eqref{fully} from~\eqref{eq: ibvp int} and using \eqref{widehatu}, we obtain 
\[
\int_{I_n}\partial_t^\alpha(\psi-\theta)\,dt
	-\int_{I_n}\mathcal{A}\theta\,dt=0.
\]
Taking the $L^2(\Omega)$-inner product with a test 
function~$\varphi\in H^1_0(\Omega)$, and applying the divergence theorem, it follows
that
\begin{equation}\label{eq: in terms of theta}
\int_{I_n} \iprod{\I^{1-\alpha} \theta',\varphi}\,dt
+\int_{I_n}\iprod{\kappa \nabla \theta,\nabla \varphi}\,dt
=\int_{I_n}\iprod{\I^{1-\alpha} \psi',\varphi}\,dt.
\end{equation}
Choosing $\varphi=\theta'$ and summing over~$n$ yields 
\[
\I(\iprod{\I^{1-\alpha}\theta',\theta'})(t_n)
	+\I(\iprod{\kappa \nabla \theta,\nabla \theta'})(t_n)
=\I(\iprod{\I^{1-\alpha} \psi',\theta'})(t_n).
\]
Since  $\I(\iprod{\kappa\nabla\theta,\nabla\theta'})(t_n)
=\tfrac12(\|\sqrt{\kappa}\nabla \theta(t_n)\|^2-\|\sqrt{\kappa}\nabla\theta(0)\|^2)=\tfrac12\|\sqrt{\kappa}\nabla \theta(t_n)\|^2$,
\begin{equation}\label{eq:varphi=theta}
\I(\iprod{\I^{1-\alpha}\theta',\theta'})(t_n)\le \I(\iprod{\I^{1-\alpha} \psi',\theta'})(t_n).
\end{equation}

To proceed in our analysis, we make use  of the following technical lemma.  
For the proof, we refer to Mustapha and Sch\"otzau~\cite[Lemma~3.1
(iii)]{MustaphaSchoetzau2014}.
\begin{lemma}\label{lem: continuity}
For $0<\alpha\le 1$ and $\epsilon >0$,
\begin{align*}
\I(\iprod{\I^{1-\alpha}v,w})(t)
	&\le \frac{1}{\alpha}\Big(
	\I(\iprod{\I^{1-\alpha}v,v})(t)\Big)^{1/2} 
	\Big(\I(\iprod{\I^{1-\alpha}w,w})(t)\Big)^{1/2}\\
&\le\epsilon\I(\iprod{\I^{1-\alpha}v,v})(t) 
	+\frac{1}{4\epsilon\alpha^2}\I(\iprod{\I^{1-\alpha}w,w})(t).
\end{align*}
\end{lemma}
For later use, by expanding  $\iprod{\mathcal{I}^{1-\alpha}(v+w),v+w)}$ then applying  Lemma~\ref{lem: continuity} with  $\epsilon=1/(2\alpha)$ we deduce the inequality in the next lemma. 
\begin{lemma}\label{lem: v plus w}
For $0<\alpha\le1$,
\[
\I(\iprod{\mathcal{I}^{1-\alpha}(v+w),v+w})(t)
    \le(1+\alpha^{-1})\Big(\I(\iprod{\mathcal{I}^{1-\alpha}v,v})(t)
    +\I(\iprod{\mathcal{I}^{1-\alpha}w,w})(t)\Big).
\]
\end{lemma}

We now apply Lemma \ref{lem: continuity} to the right-hand side 
of~\eqref{eq:varphi=theta} with~$\epsilon=1/(2\alpha^2)$. Multiplying
through by~$2$, and then cancelling the similar terms, leads to the estimate below that  will be used later in our convergence analysis.
\begin{equation}\label{theta<psi}
\I(\iprod{\I^{1-\alpha} \theta',\theta'})(t_n)
	\le\frac{1}{\alpha^2}\I(\iprod{\I^{1-\alpha}\psi',\psi'})(t_n),
\end{equation}

Under reasonable regularity assumptions, a novel error analysis involving implicit interpolations and a duality argument leads to the convergence results in the next theorem.  With $J=(0,T)$, an optimal $O(\tau^2)$-rate of
convergence is achieved in the $L^2(J; L^2(\Omega))$-norm. 
Our numerical results illustrate this in the stronger $L^\infty(J; L^2(\Omega))$-norm. Moreover, our numerical results suggest that the condition on the graded mesh exponent can be further relaxed. More precisely, instead of  $\gamma>\max\{2/\sigma, (3-\alpha)/(2\sigma-\alpha)\}$
it suffices to impose $\gamma> 2/\sigma$.

The developed error analysis requires the following regularity property
\cite[Theorems 2.1~and 2.2]{JinLiZhou2020},
\cite[Theorems 1~and 2]{SakamotoYamamoto2011}, and \cite{McLean2010}: for
some $\sigma>0$,
\begin{equation}\label{time regularity}
t\|u'(t)\|+t^2\|u''(t)\|+t^3\|u'''(t)\|\le C t^\sigma
    \quad\text{for $t>0$.}
\end{equation}
For example, if  $f\equiv 0$ and $u_0\in\dot H^r(\Omega)$ with $1\le r\le 2$, then \eqref{time regularity}
holds true for $\sigma=r\alpha/2$.

For a given time interval $Q$, the norms in $L^\infty\bigl(Q;L^2(\Omega)\bigr)$ and $L^2\bigl(Q;L^2(\Omega)\bigr)$ are  respectively denoted by $\|w\|_Q=\sup_{t\in Q}\|w(t)\|$ and $\|w\|_{L^2(Q)}=\bigl(\int_Q\|w(t)\|^2\,dt\bigr)^{1/2}$.

\begin{theorem}\label{time-convergence}
Let $u$ and $ U$ be the solutions of \eqref{eq:L1}~and \eqref{fully},
 respectively. If the graded time mesh exponent
$\gamma>\max\{2/\sigma, (3-\alpha)/(2\sigma-\alpha)\}$ and if the regularity assumption \eqref{time regularity} holds true with~$\sigma>\alpha/2$, then for $0<\alpha<1,$ we have $\|u- U\|_{L^2(J)}\le C\alpha^{-2}\tau^2$.  
\end{theorem}
{\em Proof.}
The desired estimate follows from Lemma~\ref{I1-alphapsi'psi'-estimate}~and
Theorem~\ref{etaboundpre} below. 
$\quad \Box$

\section{Errors from implicit interpolations}\label{sec:average}
In preparation for our convergence analysis, we now study the error from
approximating $u$ by~$\widehat u$, and proceed to estimate $\|\psi\|$~and
$\I(\iprod{\I^{1-\alpha} \psi',\psi'})$. These estimates assume that the
regularity property~\eqref{time regularity} holds.  For ease of reference, we
here introduce the parameter $\delta=\sigma-\frac{2}{\gamma}$
which will subsequently appear repeatedly. We start this section with  the following
representation of the implicit interpolation error in the approximation $u\approx\widehat u$ at~$t=t_n$.
\begin{lemma}\label{lem:psi-representation}
For $1\le n\le N$, $\psi^n=\sum_{j=1}^n(-1)^{n+j+1}\Delta_j$ where 
\[\Delta_j=\frac{2}{\tau_j}\int_{I_j}(u-u_I)\,dt
    =-\frac{1}{\tau_j}\int_{I_j}(t_j-t)(t-t_{j-1})u''(t)\,dt.
\]
\end{lemma}
{\em Proof.}
Since $\int_{I_n}u_I\,dt=\tfrac12\tau_n(u^n+u^{n-1})$, we find using \eqref{widehatu} that
\[
\psi^n=-\psi^{n-1}-\frac{2}{\tau_n}\int_{I_n}(u-u_I)\,dt.
\]
The formula for~$\psi^n$ then follows by induction on~$n$, after noting that $\psi^0=0$.
Recalling the Peano kernel for the trapezoidal rule, we see that
\[
\int_{I_j}(u-u_I)\,dt=-\frac{1}{2}\int_{I_j}(t_j-s)(s-t_{j-1})u''(s)\,ds,
\]
implying the second expression for~$\Delta_j$.
$\quad \Box$

\begin{lemma}\label{estimate of psin0}
For $n\ge 1$,
\[
\|\psi^n\|\le\int_{I_1}t\|u''(t)\|\,dt+\frac{1}{12}\biggl(\tau_2^2\|u''(t_1)\|
    +2\tau_n^2\|u''(t_n)\|
    +3\sum_{j=2}^n\tau_j^2\int_{I_j}\|u'''\|\,dt \biggr)
\]
and $\|\psi(t)\|\le\|u(t)-u_I(t)\|+\max\bigl(\|\psi^n\|, \|\psi^{n-1}\|\bigr)$ for $t\in I_n$.
\end{lemma}
{\em Proof.}
For $t \in I_j$, we have the identity
\[-u''(t)=\frac12\int_t^{t_j}u'''(s)\,ds-\frac12\int_{t_{j-1}}^tu'''(s)\,ds-    (u''(t_j)+u''(t_{j-1}))/2.\]
Multiply both sides by $(t_j-t)(t-t_{j-1})$ and integrate to obtain
\[\Delta_j=-\frac{\tau_j^2}{12}[u''(t_j)+u''(t_{j-1})]+R_j,\]
where
\[
R_j=\frac{1}{2\tau_j}\int_{I_j}(t_j-t)(t-t_{j-1})\biggl(
    \int_t^{t_j}u'''(s)\,ds-\int_{t_{j-1}}^tu'''(s)\,ds\biggr)\,dt.
\]
Thus, by Lemma~\ref{lem:psi-representation},
\[
(-1)^n\psi^n=\Delta_1+\frac{1}{12}\sum_{j=2}^n(-1)^j\tau_j^2
    [u''(t_j)+u''(t_{j-1})]-\sum_{j=2}^n(-1)^jR_j.
\]
Shifting the summation index,  so
\[
\sum_{j=2}^n(-1)^j\tau_j^2u''(t_{j-1})
    =\tau_2^2u''(t_1)
    -\sum_{j=2}^{n-1}(-1)^j\tau_{j+1}^2u''(t_j).\]
Since $\Delta_1=-\tau_1^{-1}\int_{I_1}(t_1-s)su''(s)\,ds$ and since $\|R_j\|\le\frac{\tau_j^2}{6}\int_{I_j}\|u'''(t)\|\,dt,$ 
\begin{align*}
\|\psi^n\|&\le\int_{I_1}t\|u''(t)\|\,dt+\frac{1}{12}\biggl(
\tau_2^2\|u''(t_1)\|+\tau_n^2\|u''(t_n)\|\\
    &\qquad{}+\sum_{j=2}^{n-1}(\tau_{j+1}^2-\tau_j^2)\|u''(t_j)\|
    +2\sum_{j=2}^n\tau_j^2\int_{I_j}\|u'''(t)\|\,dt\biggr).
\end{align*}
Using  
\[
    \sum_{j=2}^{n-1}\tau_{j+1}^2\|u''(t_j)\|
=\sum_{j=3}^n\tau_j^2\|u''(t_{j-1})\|
    \le\sum_{j=3}^n\tau_j^2\Big(\|u''(t_j)\|+\int_{I_j}\|u'''(t)\|\,dt\Big),
\]
and the bound for $\|\psi^n\|$ follows after canceling the common terms. 

The interpolant~$\psi_I$, defined as in~\eqref{eq:uI}, satisfies $\psi_I=u_I-\widehat u$, leading to the
representation
\begin{equation}\label{eq:psi-interpolate}
\psi=u-\widehat u=\psi_I+u-u_I,
\end{equation}
which implies the bound for~$\|\psi\|$ because
$\|\psi_I\|\le\max\bigl(\|\psi^n\|, \|\psi^{n-1}\|\bigr)$ on $I_n$.
$\quad \Box$

\begin{corollary}\label{estimate of psin} Under the regularity assumption
in~\eqref{time regularity} and for a time mesh of the
form~\eqref{eq: time mesh} with grading parameter~$\gamma \ge 1$ we have, 
for $n\ge 1$,
\[\|\psi^n\|\le\|\psi\|_{I_n}\le C\times\begin{cases}
      \tau^2\log(t_n/t_1),&\text{if $\gamma=2/\sigma$,}\\
    \tau^{\min(\gamma\sigma,2)}t_n^{\max(0,\delta)},  &\text{if $\gamma\ne 2/\sigma$.}
\end{cases}
\]
\end{corollary}
{\em Proof.}
First we show that for $\gamma \ge 1,$ 
\begin{equation}\label{eq:interpolation-error}
\|u-u_I\|_{I_n}\le C
    \tau^{\min(\gamma\sigma,2)}t_n^{\max(0,\delta)}. 
\end{equation}
 Since the interpolation
error~$u-u_I$ vanishes if $u$ is a polynomial of degree~$1$, by computing the
Peano kernel one finds that for~$t\in I_n$,
\begin{equation}\label{eq:uI-error}
u(t)-u_I(t)=-\frac{t_n-t}{\tau_n}\int_{t_{n-1}}^t(s-t_{n-1})u''(s)\,ds
    -\frac{t-t_{n-1}}{\tau_n}\int_t^{t_n}(t_n-s)u''(s)\,ds.
\end{equation}
The right-side is bounded by 
\[\int_{t_{n-1}}^t(s-t_{n-1})\|u''(s)\|\,ds
    +\int_t^{t_n}(t-t_{n-1})\|u''(s)\|\,ds
    \le\int_{I_n}(s-t_{n-1})\|u''(s)\|\,ds
\]
and so, by using  the time mesh properties in~\eqref{time mesh properties}, we get 
\[
\|u-u_I\|_{I_n}\le C \tau_n t_n^{-1} \int_{I_n} t\|u''(t)\|\,dt
    \le C \tau_n t_n^{-1}  \int_{I_n}t^{\sigma-1}\,dt\le C\tau_n^2t_n^{\sigma-2},~~n\ge 1.
\]
 Since $\tau_n \le C \tau t_n^{1-1/\gamma}$, the proof of~\eqref{eq:interpolation-error} is completed after noting that
\begin{equation}\label{eq:tau_and_t}
\tau_n^2t_n^{\sigma-2}\le C\tau^2t_n^\delta\le C\tau^2\max(t_n^\delta,t_1^\delta)
\le C\tau^{\min(\gamma\sigma,2)}t_n^{\max(0,\delta)}.
\end{equation}

Turning to the estimate for~$\psi^n$, Lemma~\ref{estimate of psin0} and \eqref{time regularity} imply that
\[
 \|\psi^n\|\le C\int_{I_1} t^{\sigma-1}\,dt
    +C\tau_2^2 t_1^{\sigma-2}+ C\tau_n^2 t_n^{\sigma-2}+
    C\sum_{j=2}^n \tau_j^2\int_{I_j}t^{\sigma-3}\,dt,
\]
for $1\le n\le N.$ Since $t_1=\tau^\gamma$ and  $\tau_2\le 2^\gamma \tau^\gamma,$ $\int_{I_1}t^{\sigma-1}\,dt+\tau_2^2t_1^{\sigma-2}\le C\tau^{\gamma\sigma},$
and we again bound $\tau_n^2t_n^{\sigma-2}$ using \eqref{eq:tau_and_t}.
For the sum over~$j$,
\[
\sum_{j=2}^n\tau_j^2\int_{I_j}t^{\sigma-3}\,dt
    \le C\sum_{j=2}^n\tau^2t_j^{2-2/\gamma}\int_{I_j}t^{\sigma-3}\,dt
    \le C\tau^2\int_{t_1}^{t_n}t^{\delta-1}\,dt,
\]
and the estimate for~$\|\psi^n\|$ follows.
$\quad \Box$

The next target is to estimate
$\I(\iprod{\I^{1-\alpha} \psi',\psi'})(t_n)$. Preceding this, we need to
bound $\|\psi'\|$ in the next lemma.

\begin{lemma}\label{estimate of psi'}
We have $\|\psi'(t)\|\le t^{\sigma-1}$ for~$t\in I_1$. Moreover,
$\|\psi'(t)\|\le C\tau^2\tau_n^{-1}t_n^\delta$ for $t\in I_n$ with $n\ge 2,$ and for~$\delta>0$.
\end{lemma}
{\em Proof.}
Differentiating \eqref{eq:psi-interpolate} and \eqref{eq:uI-error}, we see 
 for $t\in I_n$ that 
\[
\psi'(t)=\frac{1}{\tau_n}\int_{t_{n-1}}^t(s-t_{n-1})u''(s)\,ds
    -\frac{1}{\tau_n}\int_t^{t_n}(t_n-s)u''(s)\,ds+\tau_n^{-1}(\psi^n-\psi^{n-1}).
\]
Thus, if $t\in I_1$ then, by \eqref{time regularity}, Corollary~\ref{estimate of psin}, and the fact that $\psi^0=0,$ 
\[\|\psi'(t)\|
    \le C\tau^{-1}\int_0^ts^{\sigma-1}\,ds+C\tau_1^{-1}\int_t^{t_1}(t_1-s)s^{\sigma-2}\,ds
        +C\tau_1^{\sigma-1}\le Ct^{\sigma-1}\,.\]
 If $\delta>0$, $n\ge2$ and $t\in I_n$
then, recalling \eqref{eq:tau_and_t} and using again   \eqref{time regularity},
\[
\tau_n\|\psi'(t)\|
    \le C\tau_n\int_{I_n}t^{\sigma-2}\,dt+C \tau^2t_n^\delta
    \le C\tau_n^2t_n^{\sigma-2}+C \tau^2t_n^\delta\le C\tau^2t_n^\delta,\]
showing that $\|\psi'(t)\|\le C\tau^2\tau_n^{-1}t_n^\delta$.
$\quad \Box$
\begin{lemma}\label{I1-alphapsi'psi'-estimate}
Assume $\sigma>\alpha/2$. Then
$\big|\I(\iprod{\I^{1-\alpha} \psi',\psi'})(t_n)\big|
    \le C\tau^{3-\alpha}
    $ for $n\ge 1$, provided that $\gamma>\max\{2/\sigma,
(3-\alpha)/(2\sigma-\alpha)\}$.
\end{lemma}
{\em Proof.}
For $n=1$, the Cauchy-Schwarz inequality, Lemma~\ref{estimate of psi'}, and the assumption $\sigma>\alpha/2$ give
\begin{equation}\label{estimate of I-alphapsi'psi'I1}
\Big|\I(\iprod{\I^{1-\alpha} \psi',\psi'})(t_1)\Big|
\le C\int_0^{t_1} t^{\sigma-1}\int_0^t(t-s)^{-\alpha} s^{\sigma-1}\,ds\,dt 
    \le C\tau_1^{2\sigma-\alpha}.
\end{equation}
To deal with the case~$n\ge2$, we make the splitting
$\I^{1-\alpha}\psi'=T_1+T_2$ where
\[
T_1(t)=\int_0^{t_{j-1}}\omega_{1-\alpha}(t-s)\psi'(s)\,ds
\quad\text{and}\quad
T_2(t)=\int_{t_{j-1}}^t\omega_{1-\alpha}(t-s)\psi'(s)\,ds
\]
for $t\in I_j$~and $j\ge 2$. Using Lemma~\ref{estimate of psi'}~and
\eqref{time mesh properties},  we observe that
\[
\|T_2(t)\|
\le C\int_{t_{j-1}}^t\omega_{1-\alpha}(t-s)\,\frac{\tau^2s^\delta}{\tau_j}\,ds
\le C\frac{\tau^2t^\delta}{\tau_j}\omega_{2-\alpha}(t-t_{j-1})
\le C\tau^2t^\delta\tau_j^{-\alpha}.
\]
For estimating $T_1(t)$, integrate by parts recalling $\psi^0=0$,
\[
T_1(t)=\omega_{1-\alpha}(t-t_{j-1})\psi^{j-1}
+\int_0^{t_{j-1}}\omega_{-\alpha}(t-s) \psi(s)\,ds,
\]
where $\omega_{-\alpha}(t)=\omega_{1-\alpha}'(t)=-\alpha t^{-\alpha-1}/\Gamma(1-\alpha)$.
We apply Corollary~\ref{estimate of psin} to conclude that  
\[ \|T_1(t)\|\le C\tau^2 t_{j-1}^{\delta}\biggl(\omega_{1-\alpha}(t-t_{j-1})
    -\int_0^{t_{j-1}}\omega_{-\alpha}(t-s)\,ds\biggr)
    \le C\tau^2 t_j^{\delta}\omega_{1-\alpha}(t-t_{j-1}).
\]
Lemma~\ref{estimate of psi'} and above estimates for $\|T_1\|$~and $\|T_2\|$ yield
\[   \|\psi'\|_{I_j}\int_{I_j}\bigl(\|T_1\|+\|T_2\|\bigr)\,dt
    \le C(\tau^2t_j^{\delta}\tau_j^{-1})\tau^2t_j^\delta(\tau_j^{1-\alpha}) \le \tau^4t_j^{2\delta}\tau_j^{-\alpha},~~{\rm for}~~j\ge 2\,.
\]
By using this and \eqref{estimate of I-alphapsi'psi'I1},
and noting that $\gamma(2\sigma-\alpha)>3-\alpha$, we reach
\begin{align*}
\Big|\I(&\iprod{\I^{1-\alpha}\psi',\psi'})(t_n)\Big|
    \le\biggl|\int_0^{t_1}\iprod{\I^{1-\alpha}\psi',\psi'}\,dt\biggr|
    +\biggl|\int_{t_1}^{t_n}\iprod{\I^{1-\alpha}\psi',\psi'}\,dt\biggr|\\
&\le C\tau^{3-\alpha}+\sum_{j=2}^n
    \|\psi'\|_{I_j}\int_{I_j}\bigl(\|T_1\|+\|T_2\|\bigr)\,dt\le C\tau^{3-\alpha}+C\sum_{j=2}^n\tau^4t_j^{2\delta}\tau_j^{-\alpha},
\end{align*}
for $j\ge 2.$ By~\eqref{time mesh properties}, if $j\ge2$ then
$\tau^{1+\alpha}\le C\tau_j^{1+\alpha}t_j^{-(1+\alpha)(1-1/\gamma)}$ so
\[\tau^{1+\alpha}\sum_{j=2}^nt_j^{2\delta}\tau_j^{-\alpha}
\le C\sum_{j=2}^nt_j^{2\delta-(1+\alpha)(1-1/\gamma)}\tau_j
\le C\int_{t_1}^{t_n}t^{2\sigma-\alpha-(3-\alpha)/\gamma-1}\,dt
\le C,\]
and therefore the desired bound holds.
$\quad \Box$
\section{Errors from the time discretizations}\label{sec:semi}
This section is devoted to estimating the error~$\eta=u- U$ from the time
discretization in the norm of $L^2(J;L^2(\Omega))$.  To achieve an
optimal convergence rate, we employ a duality argument in addition to the usage of the time graded meshes.  By reversing the order
of integration, we find that
\[
\I(\iprod{\I^\alpha v,w})(T)=\I(\iprod{v,\J{T}{\alpha}w})(T)
\quad\text{where}\quad
(\J{T}{\alpha}w)(t)=\int_t^T\omega_\alpha(s-t)w(t)\,ds.
\]
Using $\I(\iprod{\partial_t^\alpha v,w})(T)
    =\I(\iprod{v',\J{T}{1-\alpha}w})(T)$ and integrating by parts yield
\begin{equation}\label{adjoint}
\I(\iprod{\partial_t^\alpha v,w})(T)
    =-\iprod{v(0),(\J{T}{1-\alpha}w)(0)}
    -\I(\iprod{v,(\J{T}{1-\alpha}w)'})(T),
\end{equation}
and since $\partial_t^\alpha v
=\I^{1-\alpha}v'=(\I^{1-\alpha}v)'-v(0)\omega_{1-\alpha}$,
\begin{equation}\label{RL-adjoint}
\I(\iprod{(\I^{1-\alpha}v)',w})(T)
    =-\I(\iprod{v,(\J{T}{1-\alpha}w)'})(T).
\end{equation}
We remark that Zhang et al.~\cite[Equation~(89)]{ZhangZengJiangZhang2023} have
recently exploited this dual operator~$w\mapsto-(\J{T}{1-\alpha}w)'$ in the
error analysis of a discontinuous Galerkin scheme for~\eqref{eq:L1}.

Suppose that $\varphi$ satisfies the final-value problem
\begin{equation}\label{dual problem}
-(\J{T}{1-\alpha}\varphi)'+\mathcal{A}\varphi=\eta(t)\quad\text{for $0<t<T$,}
\quad\text{with $\varphi(T)=0$,}
\end{equation}
subject to homogeneous Dirichlet boundary conditions. Let
$y(t)=\varphi(0)+\int_0^t\varphi(s)\,ds$ so that $y$ solves the initial-value
problem
\begin{equation}\label{eq:y}
y'=\varphi\quad\text{for~$0<t<T$,}\quad\text{with $y(0)=\varphi(0)$,}
\end{equation}
and with  $y_I$ denotes the continuous piecewise-linear function that
interpolates $y$ at the time levels~$t_j$, put
\begin{equation}\label{eq:Y}
Y=y-y_I.
\end{equation}

\begin{lemma}\label{lem:norm-eta}
With the notation above,  $\|\eta\|^2_{L^2(J)}\le\I(
    \iprod{Y',\I^{1-\alpha}\eta'+\mathcal{A}\eta})(T).$
\end{lemma}
{\em Proof.}
Using \eqref{adjoint}, \eqref{dual problem}, $\eta(0)=0$ and \eqref{eq:y},
\[
\|\eta\|^2_{L^2(J)}
=\I(\iprod{-(\J{T}{1-\alpha}\varphi)'+\mathcal{A}\varphi,\eta})(T)
=\I(\iprod{y',\I^{1-\alpha}\eta'+\mathcal{A}\eta})(T).
\]
At the same time, \eqref{eq: ibvp int} and  \eqref{fully} imply that
\[
\frac{1}{\tau_n}\int_{I_n}\bigl(
    \partial_t^\alpha\eta+\mathcal{A}\eta\bigr)\,dt=0
\quad\text{so}\quad
\int_{I_n}\iprod{y_I',\partial_t^\alpha\eta+\mathcal{A}\eta}\,dt=0,
\]
because $y_I'$ is constant on~$I_n$.  Since $Y'=y'-y_I'$, the inequality
follows at once.
$\quad \Box$

We will show below in Theorem~\ref{thm:Y} that the interpolation error~$Y$
satisfies
\begin{equation}\label{assumptions}
\I(\iprod{Y',\I^{1-\alpha}Y'})(T)
    +\I(\iprod{\mathcal{A}Y, \I^{\alpha}\mathcal{A}Y'})(T)
\le C\tau^{1+\alpha}\|\eta\|^2_{L^2(J)}.
\end{equation}
Assuming this fact for now, we can derive an estimate for~$\eta$ in terms
of~$\psi'$. For convenience, we use the notations: $F(g)=\bigl(\I(\iprod{\I^\alpha g',g})(T)\bigr)^{1/2}$ and $G(g)=\bigl(\I(\iprod{\I^{1-\alpha}g',g'})(T)\bigr)^{1/2}$.

\begin{theorem}\label{etaboundpre}
We have  $\alpha^2\|\eta\|^2_{L^2(J)}
    \le C\,\tau^{\alpha+1}
    \I(\iprod{\I^{1-\alpha} \psi',\psi'})(T).$
\end{theorem}
{\em Proof.}
It suffices to estimate the right-hand side of the inequality in
Lemma~\ref{lem:norm-eta}. By Lemma~\ref{lem: continuity},
\[
\alpha\I(\iprod{Y',\I^{1-\alpha}\psi'})(T)
\le G(Y)G(\psi)~~{\rm and}~~ 
\alpha\I(\iprod{Y',\I^{1-\alpha}\theta'})(T)
    \le G(Y)G(\theta)\,.
\]
After using \eqref{theta<psi} and $\eta'=\psi'-\theta'$, we conclude that
\begin{equation}\label{eq:eta-dash}
\alpha^2\I(\iprod{Y',\I^{1-\alpha}\eta'})(T)\le
2 G(Y)G(\psi).
\end{equation}

Since $Y(t_j)=0$ for $0\le j\le N$, integrating by parts, using
$\mathcal{A}Y=\I\mathcal{A}Y'=\I^{1-\alpha}(\I^\alpha \mathcal{A}Y')$, and applying Lemma \ref{lem: continuity}, 
\begin{align*}
\alpha \I(\iprod{Y',\mathcal{A}\theta})(T)
    &
=-\alpha\I(\iprod{\theta',\I^{1-\alpha}(\I^\alpha\mathcal{A}Y')})(T)
\le G(\theta)F({\mathcal A}Y).
\end{align*}
The same estimate holds with $\theta$ replaced by~$\psi$, so because
$\eta=\psi-\theta$, and using \eqref{theta<psi} again,
\begin{equation}\label{eq:Aeta}
\alpha^2 \I(\iprod{Y',\mathcal{A}\eta})(T)
\le 2F({\mathcal A}Y)G(\psi).
\end{equation}
Adding \eqref{eq:eta-dash}~and \eqref{eq:Aeta}, we see 
$\alpha^2\|\eta\|^2_{L^2(J)}\le2G(\psi)(G(Y)+F({\mathcal A}Y))$ by
Lemma~\ref{lem:norm-eta}.
Squaring both sides, we have
\[
\alpha^4\|\eta\|^4_{L^2(J)}\le 4(G(\psi))^2(G(Y)+F({\mathcal A}Y))^2
    \le 8 (G(\psi))^2((G(Y))^2+(F({\mathcal A}Y))^2).
\]
Since  $(G(Y))^2+(F({\mathcal A}Y))^2$ is just the left-hand side
of~\eqref{assumptions}, the desired inequality followed  after cancelling the common factor
$\|\eta\|^2_{L^2(J)}$.
$\quad \Box$

It remains to prove \eqref{assumptions}. We start by showing
preliminary bounds for $\big|\I(\iprod{Y',\I^{1-\alpha}Y'})(T)\big|$ and
$\big|\I(\iprod{\mathcal{A}Y,\I^{1-\alpha}\mathcal{A}Y'})(T)\big|$
in  the next two lemmas.

\begin{lemma}\label{assumptions1}
\[
\Big|\I(\iprod{Y', \I^{1-\alpha}Y'})(T)\Big|
   \le C(1-\alpha)\sum_{j=1}^{N-2}\tau_j^{-\alpha} \|Y\|_{I_j}^2
    +C\tau^{1-\alpha}\|Y'\|^2_{L^2(J)}.
\]
\end{lemma}
{\em Proof.}
For~$t\in I_j$ with $j\ge 2$, we write $\I^{1-\alpha} Y'(t)=S_1(t)+S_2(t)$ where
\[
S_1(t)=\int_{t_{j-2}}^t \omega_{1-\alpha}(t-s)Y'(s)\,ds
~\text{and}~
S_2(t)=\int_0^{t_{j-2}}\omega_{1-\alpha}(t-s)Y'(s)\,ds.
\]
Applying  the Cauchy-Schwarz inequality, integrating, changing the order
of integration, and integrating again, yields
\begin{align*}
\|S_1(t)&\|^2_{L^2(I_j)}\le \int_{I_j}
    \biggl(\int_{t_{j-2}}^t \omega_{1-\alpha}(t-s)\,ds\biggr)
    \biggl(\int_{t_{j-2}}^t \omega_{1-\alpha}(t-s)\|Y'(s)\|^2\,ds\biggr)\,dt\\
&\le C(\tau_j+\tau_{j-1})^{1-\alpha}\biggl(\int_{I_j}\int_s^{t_j}
+\int_{I_{j-1}}\int_{I_j}\biggr)\omega_{1-\alpha}(t-s)\,dt\,\|Y'(s)\|^2\,ds\\
&\le C(\tau_j+\tau_{j-1})^{2(1-\alpha)}
    \int_{t_{j-2}}^{t_j}\|Y'(s)\|^2\,ds
\le C\tau_j^{2(1-\alpha)}\int_{t_{j-2}}^{t_j}\|Y'(s)\|^2\,ds.
\end{align*}
For $t \in I_1$, let $S_1(t)=\int_0^t\omega_{1-\alpha}(t-s)Y'(s)\,ds$ and
$S_2(t)=0$. Following the steps above, we find that $\|S_1\|^2_{L^2(I_1)}
\le C\tau_1^{2(1-\alpha)}\|Y'\|^2_{L^2(I_1)}$.  Thus,
$\|S_1\|^2_{L^2(J)}
    \le C\tau^{2(1-\alpha)}\|Y'\|^2_{L^2(J)},$ and consequently
\begin{equation}\label{Y'S1-estimate}
\I(|\iprod{Y', S_1}|)(T)
\le\|Y'\|_{L^2(J)}\|S_1\|_{L^2(J)}
\le C\tau^{1-\alpha}\|Y'\|^2_{L^2(J)}.
\end{equation}

Turning to the second term~$S_2(t)$, we integrate by parts (noting that
$Y^{j-2}=0=Y^0$) and obtain
\[
\|S_2(t)\|=\biggl\|\int_0^{t_{j-2}}\omega_{-\alpha}(t-s)Y(s)\,ds\biggr\|
    \le\sum_{i=1}^{j-2}\|Y\|_{I_i}\int_{I_i}|\omega_{-\alpha}(t-s)|\,ds.
\]
Since $|\omega_{-\alpha}(t-s)|\le|\omega_{-\alpha}(t_{j-1}-s)|$ for~$t\in I_j$,
\[
\int_{I_j}|\iprod{Y',S_2}|\,dt
     \le\sum_{i=1}^{j-2}\int_{I_j}\|Y'(t)\|\|Y\|_{I_i}\,dt\int_{I_i}
        |\omega_{-\alpha}(t_{j-1}-s)|\,ds,
\]
with
\[
\int_{I_j}\|Y'(t)\|\|Y\|_{I_i}\,dt
\le \frac12\biggl(\|Y\|_{I_i}^2+\tau_j\|Y'\|^2_{L^2(I_j)}\biggr),
\]
and, remembering that $\omega_{-\alpha}(t)=-\alpha\,t^{-1}\omega_{1-\alpha}(t)$,
\[
\int_{I_i}|\omega_{-\alpha}(t_{j-1}-s)|\,ds
=\omega_{1-\alpha}(t_{j-1}-t_i)-\omega_{1-\alpha}(t_{j-1}-t_{i-1}).
\]
Since $t_j-t_i\ge t_{j-1}-t_{i-1}$, $\omega_{1-\alpha}(t_{j-1}-t_{i-1})
\ge \omega_{1-\alpha}(t_j-t_i)$,
and thus,
\begin{multline*}
\int_{I_j}|\iprod{Y',S_2}|\,dt\le\frac{1}{2}\sum_{i=1}^{j-2}\|Y\|_i^2
    [\omega_{1-\alpha}(t_{j-1}-t_i)-\omega_{1-\alpha}(t_j-t_i)]\\
    +\frac{\tau_j}{2}\|Y'\|^2_{L^2(I_j)}
        \int_0^{t_{j-2}}|\omega_{-\alpha}(t_{j-1}-s)|\,ds,
\end{multline*}
and in the second term, noting that
$1/\Gamma(1-\alpha)=(1-\alpha)/\Gamma(2-\alpha)$,
\[
\int_0^{t_{j-2}}|\omega_{-\alpha}(t_{j-1}-s)|\,ds
    =\omega_{1-\alpha}(t_{j-1}-t_{j-2})-\omega_{1-\alpha}(t_j)
    \le C(1-\alpha)\tau_{j-1}^{-\alpha}.
\]
By interchanging the order of the double sum,
\begin{align*}
&\I(|\iprod{Y',S_2}|)(T)
    =\sum_{j=2}^N\int_{I_j}|\iprod{Y',S_2}|\,dt\\
    &\le C(1-\alpha)\Big(\sum_{i=1}^{N-2}\|Y\|_{I_i}^2\sum_{j=i+2}^N
        [(t_{j-1}-t_i)^{-\alpha}-(t_{j}-t_{i})^{-\alpha}]
   +\Bigl(
    \max_{2\le j\le N}\tau_j\tau_{j-1}^{-\alpha}\Bigr)\|Y'\|^2_{L^2(J)}\Big)\\
    &\le C(1-\alpha)\biggl(\sum_{i=1}^{N-2}\tau_{i+1}^{-\alpha}\|Y\|_{I_i}^2
        +\tau^{1-\alpha}\|Y'\|^2_{L^2(J)}\biggr),
\end{align*}
which, combined with~\eqref{Y'S1-estimate}, yields the desired estimate.
$\quad \Box$

\begin{lemma}\label{assumption2}
\[
\Big|\I(\iprod{\mathcal{A}Y,\I^{\alpha}\mathcal{A}Y'})(T)\Big|
    \le C\sum_{j=1}^N\tau_j^\alpha\biggr[
    \biggl(\int_{I_j}\|\mathcal{A}Y'\|\,dt\biggr)^2
    +\|\mathcal{A}Y\|_{I_j}^2\biggr].
\]
\end{lemma}
{\em Proof.}
Integrating by parts,
\begin{equation}\label{eq:ass2-A}
\I(\iprod{\mathcal{A}Y,\I^{\alpha}\mathcal{A}Y'})(T)
    =-\int_0^T\iprod{\mathcal{A}Y',S_3}\,dt
     -\int_0^T\iprod{\mathcal{A}Y',S_4}\,dt,
\end{equation}
where we used the splitting $\I^{\alpha}{\mathcal A}Y=S_3+S_4$ with
\[
S_3(t)=\int_{t_{j-1}}^t\omega_\alpha(t-s)\mathcal{A}Y(s)\,ds
~~\text{and}~~
S_4(t)=\int_0^{t_{j-1}}\omega_\alpha(t-s)\mathcal{A}Y(s)\,ds
\]
for $t\in I_j$. Since $\|S_3(t)\|
\le\|\mathcal{A}Y\|_{I_j}\int_{t_{j-1}}^t\omega_\alpha(t-s)\,ds
\le C\tau_j^\alpha\|\mathcal{A}Y\|_{I_j}$,
\begin{equation}\label{eq:ass2-B}
\biggl|\int_0^T\iprod{\mathcal{A}Y',S_3}\,dt\biggr|
   \le C\sum_{j=1}^N\tau_j^\alpha\|\mathcal{A}Y\|_{I_j}
    \int_{I_j}\|\mathcal{A}Y'(t)\|\,dt.
\end{equation}
For the estimate involving $S_4$, we reverse the order of integration and then
integrate by parts, to obtain
\begin{align*}
\int_{I_j}\iprod{\mathcal{A}Y',S_4}\,dt
   &=-\int_0^{t_{j-1}}\biggiprod{\mathcal{A}Y(s),
    \int_{I_j}\omega_{\alpha-1}(t-s)\mathcal{A}Y(t)\,dt}\,ds,
   \end{align*}
and thus, applying the Cauchy-Schwarz inequality and using $\|\mathcal{A}Y\|_{I_i}\|\mathcal{A}Y\|_{I_j}\le \|\mathcal{A}Y\|_{I_i}^2+\|\mathcal{A}Y\|_{I_j}^2,$ we get  
\begin{multline*}
\biggl|\int_0^T\iprod{\mathcal{A}Y',S_4}\,dt\biggr|
\le\sum_{j=1}^N\sum_{i=1}^{j-1}
    \bigl(\|{\mathcal A}Y\|_{I_i}^2+\|\mathcal{A}Y\|_{I_j}^2\bigr)
    \int_{I_i}\int_{I_j}|\omega_{\alpha-1}(t-s)|\,dt\,ds\\
=\sum_{i=1}^{N-1}\|\mathcal{A}Y\|_{I_i}^2
    \int_{I_i}\int_{t_i}^{t_n}|\omega_{\alpha-1}(t-s)|\,dt\,ds
+\sum_{j=1}^N \|\mathcal{A}Y\|_{I_j}^2
    \int_{I_j}\int_0^{t_{j-1}}|\omega_{\alpha-1}(t-s)|\,ds\,dt\\
\le \sum_{i=1}^{N-1}\|{\mathcal A}Y\|_{I_i}^2
    \int_{I_i}\omega_{\alpha}(t_i-s)\,ds
    +\sum_{j=1}^N \|{\mathcal A}Y\|_{I_j}^2
    \int_{I_j}\omega_{\alpha}(t-t_{j-1})\,dt
    \le  C\sum_{j=1}^N\tau_j^\alpha\|{\mathcal A}Y\|_{I_j}^2.
\end{multline*}
The proof is concluded by inserting this  and \eqref{eq:ass2-B} into the
splitting~\eqref{eq:ass2-A}.
$\quad \Box$

By using the achieved estimates in the previous two lemmas, we are now able in
the next theorem to provide the missing part in the proof of
Theorem~\ref{etaboundpre}.

\begin{theorem}\label{thm:Y}
The inequality \eqref{assumptions} is satisfied by the function $Y$ defined via
\eqref{dual problem}, \eqref{eq:y} and \eqref{eq:Y}.
\end{theorem}
{\em Proof.}
Recall that $Y=y-y_I$ where $y_I$ is the piecewise linear polynomial that
interpolates $y$ at the time nodes, and $y'=\varphi$. Thus, if $t\in I_j$ then
\begin{equation}\label{eq:Y varphi}
Y(t)=y(t)-y_I(t)=\frac{t_j-t}{\tau_j}\int_{t_{j-1}}^t\varphi(s)\,ds
    -\frac{t-t_{j-1}}{\tau_j}\int_t^{t_j}\varphi(s)\,ds
\end{equation}
so $\|Y(t)\|\le\int_{I_j}\|\varphi\|\,ds$.  Similarly, replacing $u$ with $y$ in \eqref{eq:uI-error}, we have $\|Y(t)\|\le\tau_j\int_{I_j}\|\varphi'\|\,ds$. Therefore
\begin{equation}\label{eq:Y varphi estimate}
\|Y\|_{I_j}^2\le\tau_j\|\varphi\|^2_{L^2(I_j)}
\quad\text{and}\quad
\|Y\|_{I_j}^2\le\tau_j^3\|\varphi'\|^2_{L^2(I_j)}.
\end{equation}
Consider the linear operator~$B$ defined by
$(B\varphi)(t)=\tau_j^{-\alpha/2}\|Y\|_{I_j}$ for $t\in I_j$ and $1\le j\le N$. The estimates~\eqref{eq:Y varphi estimate} give
\[
\|B\varphi\|^2_{L^2(J)}\le C\tau^{1-\alpha}\|\varphi\|^2_{L^2(J)}
\quad\text{and}\quad
\|B\varphi\|^2_{L^2(J)}\le C\tau^{3-\alpha}\|\varphi'\|^2_{L^2(J)},
\]
and, since $\varphi(T)=0$
and $(\tau^{1-\alpha})^{1-\alpha}(\tau^{3-\alpha})^\alpha=\tau^{1+\alpha}$, we
may apply Corollary~\ref{cor:interpolation} from the Appendix to deduce that
\begin{equation}\label{eq:sum tau alpha Y}
\sum_{j=1}^N\tau_j^{-\alpha}\|Y\|_{I_j}^2
    =\|B\varphi\|^2_{L^2(J)}
    \le C\tau^{1+\alpha}\|(\J{T}{1-\alpha}\varphi)'\|^2_{L^2(J)}.
\end{equation}

Furthermore, by differentiating \eqref{eq:Y varphi} and \eqref{eq:uI-error},
\[
Y'(t)=\varphi(t)-\frac{1}{\tau_j}\int_{I_j}\varphi(s)\,ds
    =\frac{1}{\tau_j}\int_{t_{j-1}}^t(s-t_{j-1})\varphi'(s)\,ds
    -\frac{1}{\tau_j}\int_t^{t_j}(t_j-s)\varphi'(s)\,ds
\]
for~$t\in I_j$, so
\[
\|Y' (t)\|\le\|\varphi(t)\|+\frac{1}{\tau_j}\int_{I_j}\|\varphi(s)\|\,ds
\quad\text{and}\quad
\|Y'(t)\|\le\int_{I_j}\|\varphi'(s)\|\,ds,
\]
implying that
\begin{equation}\label{eq:Y' varphi}
\|Y'\|^2_{L^2(I_j)}\le4\|\varphi\|^2_{L^2(I_j)}
\quad\text{and}\quad
\|Y'\|^2_{L^2(I_j)}\le\tau_j^2\|\varphi'\|^2_{L^2(I_j)}.
\end{equation}
After summing over~$j$ and once again applying
Corollary~\ref{cor:interpolation}, we arrive at
\begin{equation}\label{eq:Y-RL}
\|Y'\|^2_{L^2(J)}
    \le C\tau^{2\alpha}\|(\J{T}{1-\alpha}\varphi)'\|^2_{L^2(J)}.
\end{equation}
Now take the inner product of~\eqref{dual problem}
with~$-(\J{T}{1-\alpha}\varphi)'$ in~$L_2(\Omega)$, and then integrate in time
to obtain
\begin{align*}
\|(\J{T}{1-\alpha}\varphi)'\|^2_{L^2(J)}
    -\I(\iprod{\mathcal{A}\varphi,(\J{T}{1-\alpha}\varphi)'})(T)
    &=-\I(\iprod{\eta,(\J{T}{1-\alpha}\varphi)'})(T)\\
    &\le\frac{1}{2}\|\eta\|^2_{L^2(J)}
    +\frac{1}{2}\|(\J{T}{1-\alpha}\varphi)'\|^2_{L^2(J)}.
\end{align*}
By~\eqref{RL-adjoint},
\begin{equation}\label{eq:positivity}
-\I(\iprod{\mathcal{A}\varphi,(\J{T}{1-\alpha}\varphi)'})(T)
=\I(\iprod{(\I^{1-\alpha}\mathcal{A}^{1/2}\varphi)',
\mathcal{A}^{1/2}\varphi})(T)\ge0
\end{equation}
and therefore
$\|(\J{T}{1-\alpha}\varphi)'\|^2_{L^2(J)}\le\|\eta\|^2_{L^2(J)}.$
Combining this with \eqref{eq:sum tau alpha Y}, \eqref{eq:Y-RL}, we conclude that
\[
\sum_{j=1}^N\tau_j^{-\alpha}\|Y\|_{I_j}^2
    \le C\tau^{1+\alpha}\|\eta\|^2_{L^2(J)}
\quad\text{and}\quad
\|Y'\|^2_{L^2(J)}\le C\tau^{2\alpha}\|\eta\|^2_{L^2(J)}
\]
and so, applying Lemma~\ref{assumptions1},
\begin{equation}\label{eq:first term}
\I(\iprod{Y',\I^{1-\alpha}Y'})(T)
    \le C\tau^{1+\alpha}\|\eta\|^2_{L^2(J)}.
\end{equation}
By taking the inner product of~\eqref{dual problem} with~$\mathcal{A}\varphi$ and proceeding as above, we deduce that 
\begin{equation}\label{eq:A-varphi-eta}
\|\mathcal{A}\varphi\|^2_{L^2(J)}\le\|\eta\|^2_{L^2(J)}.
\end{equation}

Repeating the arguments leading to the first estimate
in~\eqref{eq:Y' varphi} but with $Y'$ replaced by $\mathcal{A}Y'$, we see that
$\|\mathcal{A}Y'\|^2_{L^2(J)}\le \|\mathcal{A}\varphi\|^2_{L^2(J)}$ and so
\[
\sum_{j=1}^N\tau_j^\alpha\biggl(\int_{I_j}\|\mathcal{A}Y'\|\,dt\biggr)^2
    \le\sum_{j=1}^N\tau_j^{1+\alpha}\|\mathcal{A}Y'\|^2_{L^2(I_j)}
    \le C\tau^{1+\alpha}\|\mathcal{A}\varphi\|^2_{L^2(J)}.
\]
Likewise, $\|\mathcal{A}Y\|_{I_j}^2
\le\tau_j\|\mathcal{A}\varphi\|^2_{L^2(I_j)}$ by the arguments leading
to~\eqref{eq:Y varphi}, so
\[
\sum_{j=1}^N\tau_j^\alpha\|\mathcal{A}Y\|_{I_j}^2
    \le\sum_{j=1}^N\tau_j^{1+\alpha}\|\mathcal{A}\varphi\|^2_{L^2(I_j)}
    \le C\tau^{1+\alpha}\|\mathcal{A}\varphi\|^2_{L^2(J)}.
\]
Hence, by Lemma~\ref{assumption2} and \eqref{eq:A-varphi-eta},
\begin{equation}\label{eq:second term}
\I(\iprod{\mathcal{A}Y,\I^\alpha\mathcal{A}Y'})(T)
    \le C\tau^{1+\alpha}\|\eta\|^2_{L^2(J)}.
\end{equation}
Together, \eqref{eq:first term}~and \eqref{eq:second term} imply the desired
estimate~\eqref{assumptions}.
$\quad \Box$


\section{A fully discrete scheme and error analysis}\label{sec:fully}
In this section, we discretize the time-stepping scheme~\eqref{fully} in
space using the continuous piecewise-linear Galerkin FEM and hence define
a fully-discrete method. Thus, we introduce a family of regular (conforming)
triangulations~$\mathcal{T}_h$ of the domain~$\overline{\Omega}$ indexed by
$h=\max_{K\in \mathcal{T}_h}(h_K)$, where $h_{K}$ denotes the diameter of the
element~$K$. Let $V_h$ denote the  space of continuous, piecewise-linear
functions with respect to~$\mathcal{T}_h$ that vanish on~$\partial\Omega$. Let $\W(V_h)\subset {\mathcal C}([0,T];V_h)$ denote the space of linear polynomials on~$\overline I_n$ for~$1\le  n\le N$, with coefficients in~$V_h$. Motivated by the weak formulation of~\eqref{fully}, our fully-discrete
solution~$ U_h\in\W(V_h)$ is defined by requiring
\begin{equation} \label{CNFEM}
\biggiprod{\int_{I_n}\partial_t^\alpha  U_h\,dt,v_h}
    +\tau_n\iprod{\kappa\nabla  U_h^{n-1/2},\nabla v_h}
=\tau_n\iprod{\bar f_n,v_h}\quad\text{for all $v_h\in V_h$,}
\end{equation}
and for $1\le n\le N$, where $ U_h^n:= U_h(t_n)$ and
$ U_h^{n-1/2}=\tfrac12( U_h^n+ U_h^{n-1})$.  For the discrete
initial data, we choose $ U_h^0=R_h u_0\in V_h$, where $R_h:H^1_0(\Omega) \mapsto V_h$ is the Ritz projection defined by $\iprod{\kappa\nabla(R_h w-w),\nabla v_h}= 0$ for all $v_h\in V_h$.

In the next theorem, we prove that the numerical solution defined
by~\eqref{CNFEM} is second order accurate in both time and space, provided that
the time mesh exponent~$\gamma$ is chosen appropriately. In comparison, under
heavier regularity assumptions and stronger graded meshes, convergence of
order~$h^2+\tau^{(3-\alpha)/2}$ was proved by Mustapha et
al.~\cite{MustaphaAbdallahFurati2014}. Furthermore, the proof therein is more
technical and lengthy.  Use of the piecewise-linear polynomial
function~$\widehat u$, see~\eqref{widehatu}, and a duality argument allowed us to improve the convergence rate, simplify the proof and also relax the regularity assumptions. In addition to
the regularity assumption in~\eqref{time regularity}, for $t>0,$ we impose
\begin{equation}\label{spatial regularity}
\|u'(t)\|_2+t\|u''(t)\|_2\le C\,t^{\upsilon-1},\quad\text{for some~$\upsilon>0.$}
\end{equation}

\begin{theorem}\label{Convergence theorem time-space}
Let $u$ be the solution of~\eqref{eq:L1} and let $ U_h^n$ be the approximate
solution defined by~\eqref{CNFEM}.  Assume that the regularity assumptions in \eqref{time regularity}~and \eqref{spatial regularity} are satisfied for~$\sigma,\upsilon>\alpha/2$, and choose the mesh grading
exponent~$\gamma>\max\{2/\sigma,1/\upsilon,  (3-\alpha)/(2\sigma-\alpha)\}$.
Then, 
\[
\|u- U_h\|_{L^2(J)}\le C(\tau^2+h^2).
\]
\end{theorem}
{\em Proof.}
Decompose the error as $u- U_h= (u- u_h)+(u_h- U_h)$, where $ u_h$ is the Galerkin finite element solution of \eqref{eq:L1} defined by
\begin{equation}\label{FEM}
    \iprod{\partial_t^\alpha  u_h,v_h}
    +\iprod{\kappa\nabla  u_h,\nabla v_h}
=\iprod{f,v_h}\quad\text{for all $v_h\in V_h$,}
\end{equation}
for each fixed $t >0,$ with $ u_h(0)= U_h^0=R_hu_0.$  From this, the weak formulation of \eqref{eq:L1}, and the orthogonality property of the Ritz projection, we have 
\[ \iprod{\partial_t^\alpha (u_h-R_h u),v_h}
    +\iprod{\kappa\nabla  (u_h-R_hu),\nabla v_h}
=\iprod{\partial_t^\alpha (u-R_h u),v_h}\quad\text{for $v_h\in V_h$,}\]
Choose $v_h= (u_h-R_hu)'$, integrate over $(0,t)$ and apply Lemma~\ref{lem: continuity} to the right-hand side with $\epsilon =1/(4\alpha^2)$. After canceling the common terms,  
\[4\alpha^2\|\sqrt{\kappa} \nabla(u_h-R_h u)(t)\|^2
    \le \I( \iprod{\partial_t^\alpha e_h,e_h'})(t),~~{\rm with}~~ e_h=u-R_h u\,. \]
The error bound for the Ritz projection and the regularity assumption in \eqref{spatial regularity} yield  
$\|e_h'(t)\|\le Ch^2\|u'(t)\|_2\le  Ch^2 t^{\upsilon-1}$. Hence,    
$\|\partial_t^\alpha e_h(t)\|
    \le Ch^2 t^{\upsilon-\alpha}$ and consequently, $\I( \iprod{\partial_t^\alpha e_h,e_h'})(t)\le Ch^4$ for $\upsilon>\alpha/2$. Inserting this estimate in the above equation, we obtain 
$\|\nabla(u_h-R_h u)(t)\| \le C h^2$
for $\upsilon>\alpha/2$, and  thus, by applying the  Poincar\'e and triangle  inequalities, we get  
 $\|u(t)- u_h(t)\|\le Ch^2.$ 
 
 The remaining target now is to estimate   $ U_h- u_h.$  By analogy with our earlier splitting~\eqref{eq:eta-decomposition}, we let 
 \[\eta_h= u_h- U_h,\quad \psi_h=u-R_h\widehat u\quad{\rm and}\quad 
\theta_h= U_h-R_h \widehat u.\]
From \eqref{CNFEM} and \eqref{FEM}, and with $\chi_h=u_h-R_h\widehat u,$ we have 
\[\int_{I_n}[\iprod{\partial_t^\alpha \theta_h(t),v_h}
+\iprod{\kappa\nabla\theta_h(t),\nabla v_h}]\,dt
=\int_{I_n}[\iprod{\partial_t^\alpha\chi_h(t),v_h}
    +\iprod{\kappa\nabla\chi_h(t),\nabla v_h}]\,dt,\]
for $v_h\in V_h$. Using the orthogonality property of the Ritz projection and the definition
of~$\widehat u$, in addition to \eqref{FEM}, 
\[\int_{I_n}\iprod{\kappa\nabla\chi_h(t),\nabla v_h}\,dt
    =\int_{I_n}\iprod{\kappa\nabla( u_h- u)(t),\nabla v_h}\,dt=\int_{I_n}\iprod{\partial_t^\alpha(u- u_h)(t),v_h}\,dt,\]
and hence
\[\int_{I_n}[\iprod{\partial_t^\alpha\theta_h(t),v_h}
    +\iprod{\kappa\nabla\theta_h(t),\nabla v_h}]\,dt
=\int_{I_n}\iprod{\partial_t^\alpha\psi_h(t),v_h}\,dt,
\quad\text{for $v_h\in V_h$.}\]
By repeating the steps from~\eqref{eq: in terms of theta} to~\eqref{theta<psi},
we deduce that
\begin{equation}\label{eq:theta_h psi_h}
\alpha^2\I(\iprod{\partial_t^\alpha\theta_h,\theta_h'})(T)
    \le\I(\iprod{\partial_t^\alpha\psi_h,\psi_h'})(T).
\end{equation}
Applying Lemma~\ref{lem: v plus w} with $v=\widehat e_h':=(\widehat u-R_h\widehat u)'$~and
$w=\psi'$ so that $\psi_h'=v+w$, 
\[\I(\iprod{\partial_t^\alpha\psi_h,\psi_h'})(T)
    \le(1+\alpha^{-1})\Big(
    \I(\iprod{\partial_t^\alpha \widehat e_h,
    \widehat e_h'})(T)
    +\I(\iprod{\partial_t^\alpha\psi,\psi'})(T)\Big).\]
For $t\in I_j$, since $\widehat u'(t)=\tau_j^{-1}(\widehat u^j-\widehat u^{j-1}),$ the Ritz projection error bound gives
\[\|\widehat e_h'(t)\|\le Ch^2\|\widehat u'(t)\|_2
\le Ch^2\tau_j^{-1}\|\psi^j-\psi^{j-1}\|_2
    + Ch^2\tau_j^{-1}\int_{I_j}\|u'(s)\|_2\,ds\,.\]
From Lemma \ref{lem:psi-representation},  \eqref{spatial regularity} and the time mesh property \eqref{time mesh properties}, we have 
\begin{multline*}
\|\psi^j-\psi^{j-1}\|_2\le \|\psi^j\|_2+\|\psi^{j-1}\|_2\le 2\sum_{i=1}^j \int_{I_j}(t-t_{j-1})\|u''(t)\|_2\,dt
\\
\le   C\tau_1^\upsilon + C\tau \int_{t_1}^{t_j} t^{\upsilon-1/\gamma-1}\,dt
\le   C\tau_1^\upsilon + C\tau t_j^{\upsilon-1/\gamma}\le C\tau t_j^{\upsilon-1/\gamma}
\le C\tau_j t_j^{\upsilon-1},
\end{multline*}
for $\gamma>1/\upsilon$  with $j\ge 1.$ Combining the above two estimates, we conclude that 
$\|\widehat e_h'(t)\|\le Ch^2 t_j^{\upsilon-1}\le Ch^2 t^{\upsilon-1}$ for $t \in I_j$. This leads to  
$\|\widehat e_h'(t)\|\le Ch^2\omega_\upsilon(t)$ and 
$\|\partial_t^\alpha\widehat e_h(t)\|
    \le Ch^2\omega_{1-\alpha+\upsilon}(t).$
Thence,
\[
\I(\iprod{\partial_t^\alpha \widehat e_h,
    \widehat e_h '})(T)
    \le Ch^4\int_0^T\omega_{1-\alpha+\upsilon}\,\omega_\upsilon\,dt
    \le Ch^4\int_0^T t^{2\upsilon-\alpha-1}\,dt\le Ch^4,
\]
 for~$\upsilon>\alpha/2$, so by Lemma~\ref{I1-alphapsi'psi'-estimate},
\begin{equation}\label{eq:estimatepsih}
\I(\iprod{\partial_t^\alpha\psi_h,\psi_h'})(T)
    \le C\Big(h^4+\I(\iprod{\partial_t^\alpha\psi,\psi'})(T)
    \Big)
    \le C(h^4+\tau^{3-\alpha}).    
\end{equation}
Adapting \eqref{dual problem}, suppose that $\varphi_h(t) \in V_h$ satisfies the discrete final-value problem
\begin{equation}\label{dual problem discrete}
-(\J{T}{1-\alpha}\varphi_h)'(t)+\mathcal{A}_h \varphi_h(t) =\eta_h(t) \quad\text{for $0<t<T$,}
\quad\text{with $\varphi_h(T)=0$,}
\end{equation}
where the discrete elliptic operator ${\mathcal A}_h :V_h \to V_h$ is defined by $\iprod{{\mathcal A}_h v_h, q_h}=\iprod{\kappa \nabla v_h,\nabla q_h}$  for all $v_h,q_h \in V_h$. We now repeat the step in the error analysis of Section~\ref{sec:semi}, with $\eta_h,$ $\theta_h$, $\psi_h$ and $\mathcal{A}_h$ playing the roles of $\eta,$ $\theta$ $\psi$ and $\mathcal{A}$, respectively, and using \eqref{CNFEM}~and \eqref{eq:theta_h psi_h} instead of \eqref{fully}~and \eqref{theta<psi}, and  \eqref{eq:estimatepsih} instead of Lemma~\ref{I1-alphapsi'psi'-estimate}. We notice for  
$\gamma> \max\{2/\sigma,1/\upsilon, (3-\alpha)/(2\sigma-\alpha)\}$ and for $\sigma,\upsilon > \alpha/2,$ that 
\[
\| u_h- U_h\|_{L^2(J)}^2 \le C\tau^{1+\alpha}(h^4+\tau^{3-\alpha}).
\]
The proof of this theorem is completed. 
$\quad \Box$
\section{Numerical results}\label{Sec: Numeric}

In this section, we illustrate numerically the theoretical finding in
Theorem~\ref{time-convergence}.  An $O(h^2)$ convergence of the
finite element solution was confirmed for various choices of the given
data~\cite{JinZhou2017,JinLiZhou2020,Mustapha2018}. In time, some
numerical convergence results (piecewise linear discontinuous Petrov-Galerkin
method) were also delivered~\cite{MustaphaAbdallahFurati2014}. However, we
illustrate the errors and convergence rates in the stronger
$L^\infty((0,T),L^2(\Omega))$-norm on more realistic examples. We choose
$\kappa=1$,  $\Omega=(0,1)$ and a uniform spatial grid~$\mathcal{T}_h$.
In both examples, we choose $h$ so that the error from the time discretization
dominates.
  
\paragraph{Example 1.} We choose $u_0(x)=x(1-x)$~and $f\equiv0$. Thus, by
separating variables, the continuous solution has a series representation
in terms of the Mittag-Leffler function~$E_\alpha$,
\[
u(x,t)=8\sum_{m=0}^\infty\lambda_m^{-3}E_{\alpha}(-\lambda_m^2 t^{\alpha})
\sin(\lambda_m x),\quad
\text{where $\lambda_m=(2m+1)\pi$.}
\]
Since $u_0\in \dot H^{2.5^-}(\Omega)$, the regularity
estimate~\eqref{time regularity} is satisfied for~$\sigma=\alpha$.  Thus, we
expect from Theorem~\ref{time-convergence} that
$e_\tau:=\|u- U_h\|_{L^2(J)}\le C\tau^2$ provided that the mesh exponent
$\gamma> \max\{2/\alpha,(3-\alpha)/\alpha\}=(3-\alpha)/\alpha$. The
numerical results in Table~\ref{Tab 1} indicate order~$\tau^2$ convergence
in the stronger $L^\infty((0,T),L^2(\Omega))$-norm for~$\gamma> 2/\alpha$.  Rates of order~$\tau^{\sigma\gamma}$ for $1\le \gamma \le 2/\sigma=2/\alpha$ are
observed.  Thence, our imposed assumption on $\gamma$ is not sharp.

 To measure
the $L^\infty((0,T),L^2(\Omega))$~error $E_\tau:=\max_{0\le t\le T}\|u-U_h\|$, we approximated $E_\tau$
by $\max_{1\le j\le N}\max_{1\le i\le 3} \|u(t_{i,j})- U_h(t_{i,j})\|$ where
$t_{i,j}:=t_{j-1}+i\tau_j/3$. In our calculations, the $L^2(\Omega)$ norm, $\|\cdot\|,$ is approximated using the
two-point composite Gauss quadrature rule. 

\begin{figure}
\begin{center}
\includegraphics[width=7cm, height=6cm]{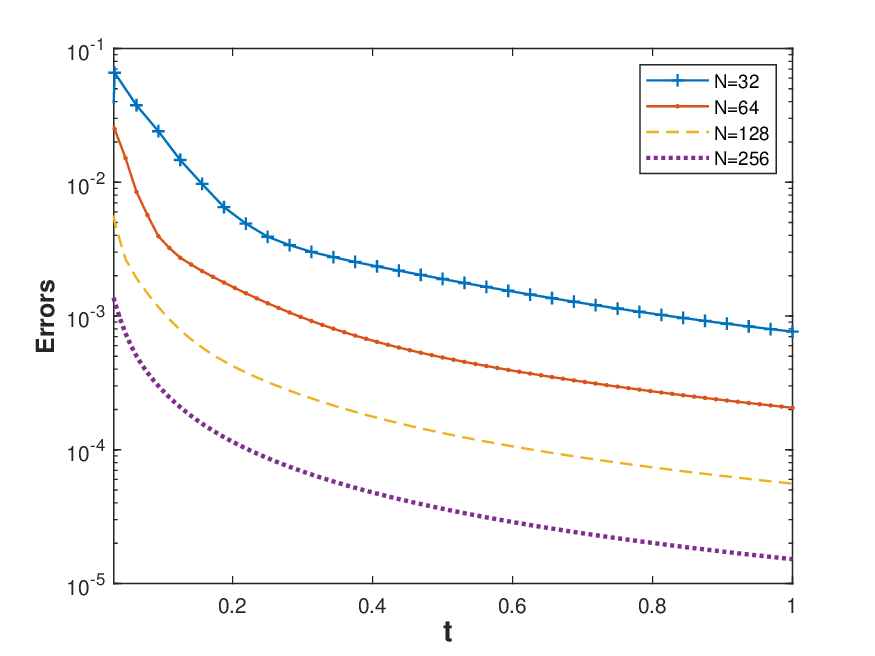}\includegraphics[width=7cm, height=6cm]{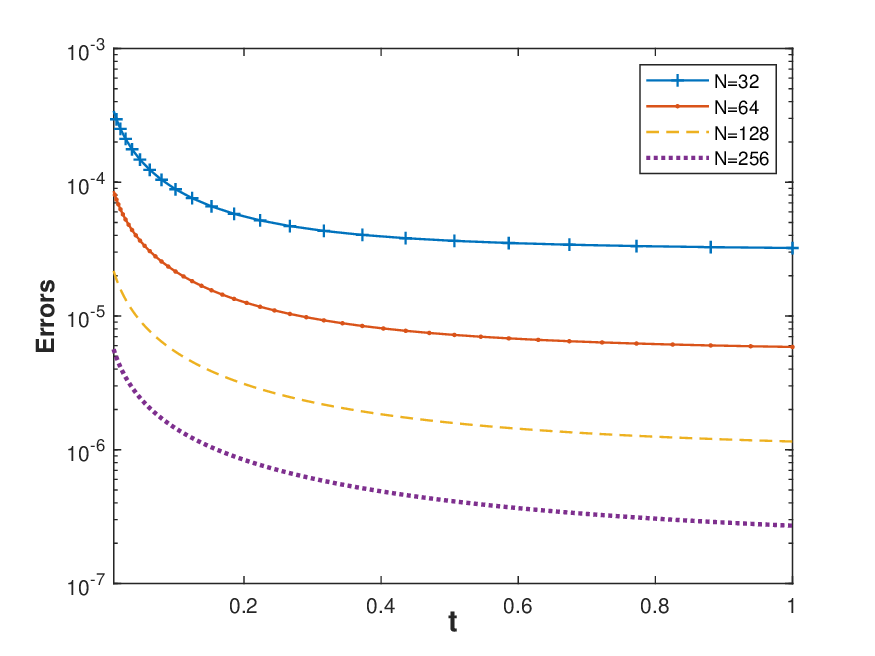}
\caption{Errors for Example~1 as functions of~$t$ for different choices of~$N$
when $\alpha=0.5$, taking $\gamma=1$ in the left figure and $\gamma=4$ in the
right figure.}
\label{Fig 1}
\end{center}
\end{figure}

In all tables and figures, we evaluated the series solution by truncating the
infinite series after $60$~terms.  The empirical convergence rate~CR is
calculated by halving~$\tau$, that is, $\text{CR}=\log_2(E_{\tau}/E_{\tau/2})$.
Figure~\ref{Fig 1} plots the nodal errors~$\|U^n_h-u(t_n)\|$
against~$t_n\in[0,1]$ for different values of~$N$ in the cases
$\gamma=1$~and $\gamma=4$.  The practical benefit of the mesh grading is
evident.

\begin{table}[ht]
\begin{center}
\begin{tabular}{|c|cc|cc|cc|cc|}
\hline 
 & \multicolumn{2}{c|}{$\gamma=1$}
& \multicolumn{2}{c|}{$\gamma=2$}
& \multicolumn{2}{c|}{$\gamma=3$}
& \multicolumn{2}{c|}{$\gamma=4$}\\
\hline
$N$&\multicolumn{2}{c|}{$E_\tau$~~~~~CR }
&\multicolumn{2}{c|}{$E_\tau$~~~~~CR}
&\multicolumn{2}{c|}{$E_\tau$~~~~~CR}
&\multicolumn{2}{c|}{$E_\tau$~~~~~CR}\\
\hline
  8& 1.011e-01&        &   5.001e-02&        &   1.861e-02&        &   7.478e-03&        \\
 16& 8.337e-02&   0.279&   2.633e-02&   0.926&   6.464e-03&   1.526&   2.090e-03&   1.839\\
 32& 6.588e-02&   0.340&   1.309e-02&   1.008&   2.259e-03&   1.517&   5.497e-04&   1.927\\
 64& 5.001e-02&   0.397&   6.464e-03&   1.018&   7.932e-04&   1.501&   1.449e-04&   1.923\\
128& 3.672e-02&   0.446&   3.205e-03&   1.012&   2.817e-04&   1.494&   3.801e-05&   1.931\\
      \hline
\end{tabular}     
\caption{Errors and empirical convergence rates for Example~1
with~$\alpha=0.5$, using different choices of the time mesh-grading
exponent~$\gamma$.} \label{Tab 1}
\end{center}
\end{table}

\begin{figure}[ht]
\begin{center}
\includegraphics[width=8cm, height=6cm]{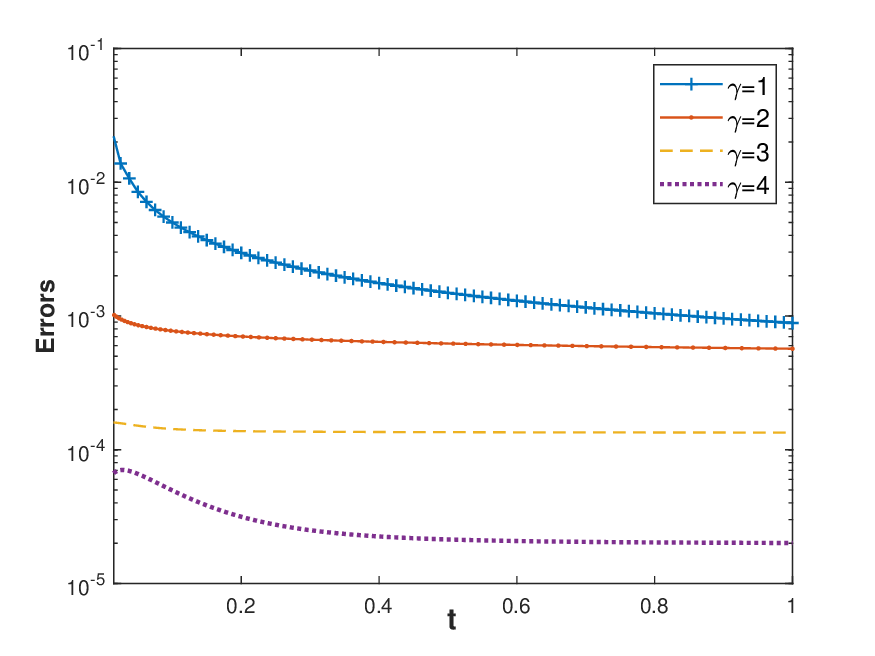}
\caption{Error at $t_n$ for $1\le n\le N$ in Example~2, for a fixed $N=80$  and different choices of the mesh exponent $\gamma$ with $\alpha=0.7$.}
\label{Fig 2}
\end{center}
\end{figure}

\paragraph{Example~2.}
We again take $f\equiv0$ but now choose less regular initial data, namely,
the hat function on the unit interval, $u_0(x)=1-2|x-\tfrac12|$. So, 
\[
u(x,t)=4\sum_{m=0}^\infty (-1)^m\lambda_m^{-2}E_\alpha(-\lambda_m^2t^{\alpha})
\sin(\lambda_mx).
\]
Since $u_0 \in \dot H^{1.5^-}(\Omega)$, the regularity
property~\eqref{time regularity} is satisfied for~$\sigma=\tfrac34\alpha$.
As in  Example 1, the numerical results in Table~\ref{Tab 2}
exhibit convergence of order~$\tau^{\sigma \gamma}$
for~$1\le\gamma\le 2/\sigma$ in the stronger $\|\cdot\|_J$-norm.  For a
graphical illustration of the impact of the graded mesh on the pointwise error,
we fixed $N=80$ in Figure \ref{Fig 2} and plotted the error at the time nodal
points for different choices of~$\gamma$.
\begin{table}[ht]
\begin{center}
\begin{tabular}{|c|cc|cc|cc|cc|}
\hline 
 & \multicolumn{2}{c|}{$\gamma=1$}
& \multicolumn{2}{c|}{$\gamma=2$}
& \multicolumn{2}{c|}{$\gamma=3$}
& \multicolumn{2}{c|}{$\gamma=4$}\\
\hline
$N$&\multicolumn{2}{c|}{$E_\tau$~~~~~CR }
&\multicolumn{2}{c|}{$E_\tau$~~~~~CR}
&\multicolumn{2}{c|}{$E_\tau$~~~~~CR}
&\multicolumn{2}{c|}{$E_\tau$~~~~~CR}\\
\hline
  8& 9.958e-02&        &   2.771e-02&        &   9.365e-03&        &   6.670e-03&        \\
 16& 6.501e-02&   0.615&   1.345e-02&   1.042&   3.144e-03&   1.575&   1.710e-03&   1.963\\
 32& 4.158e-02&   0.644&   6.509e-03&   1.048&   1.056e-03&   1.574&   4.300e-04&   1.992\\
 64& 2.771e-02&   0.586&   3.144e-03&   1.050&   3.555e-04&   1.571&   1.095e-04&   1.973\\
128& 1.924e-02&   0.526&   1.519e-03&   1.049&   1.193e-04&   1.576&            &        \\
            \hline            
\end{tabular}     
\caption{Errors and empirical convergence rates for Example~2 using different values of $N$ and different choices of the time mesh exponent $\gamma,$ with $\alpha=0.7.$}
\label{Tab 2}
\end{center}
\end{table}

\section{Appendix: An $\alpha$-robust interpolation estimate}
The purpose of this appendix is to prove Corollary~\ref{cor:interpolation}, which was used in the proof of Theorem~\ref{thm:Y}. Suppose that $Y$ is a complex Hilbert space with inner product~$(\cdot,\cdot)_Y$ and norm~$\|\cdot\|_Y$. Put
\begin{equation}\label{eq: X0 X1}
X_0=L_2(0,T;Y)\quad\text{and}\quad
X_1=\{\,u\in X_0:\text{$u'\in X_0$ and $u(0)=0$}\,\},
\end{equation}
equipped with the inner products
\[
(u,v)_{X_0}=\int_0^T\bigl(u(t),v(t)\bigr)_Y\,dt
\quad\text{and}\quad
(u,v)_{X_1}=\int_0^T\bigl(u'(t),v'(t)\bigr)_Y\,dt.
\]
We define an unbounded operator~$D_0$ on~$X_0$ with domain~$X_1$, by setting
$D_0u=u'$ on the interval~$(0,T)$.  Let $\Re z$ denote the real part of a
complex number~$z$. Since
$\Re\bigl(u(t),D_0u(t)\bigr)_Y=\tfrac12(d/dt)\|u(t)\|_Y^2$, we have
$\Re(u,D_0u)_{X_0}=\tfrac12\|u(T)\|_Y^2\ge0$ for all $u\in X_1$. In addition, one easily verifies that
\begin{equation}\label{eq: inverse lambda plus D}
(\lambda+D_0)^{-1}u(t)=\int_0^te^{-\lambda(t-s)}u(s)\,ds
\quad\text{for $0\le t\le T$ and $\lambda\in\mathbb{C}$,}
\end{equation}
showing that $D_0$ is m-accretive~\cite[p.~9]{CarracedoAlix2001}, and therefore a positive
operator~\cite[pp.~1, 11]{CarracedoAlix2001} with
$\sup_{\lambda>0}\|\lambda(\lambda+D_0)^{-1}\|_{X_0\to X_0}=1$, where $\|\cdot\|_{X_0\to X_0}$ denotes the operator norm induced by the norm of~$X_0$.

Suppose that $0<\Re\beta<1$. The operator $D_0^{-\beta}$ is defined via the integral~\cite[Equation~(1.3)]{Ashyralyev2009}
\[
D_0^{-\beta}u=\frac{\sin\pi\beta}{\pi}\int_0^\infty\lambda^{-\beta}
    (\lambda+D_0)^{-1}u\,d\lambda
\quad\text{for $u\in X_0$.}
\]
A short calculation using~\eqref{eq: inverse lambda plus D} then 
shows~\cite[Theorem~2.2]{Ashyralyev2009} that
$D_0^{-\beta}u=\mathcal{I}^\beta u$ for $0<\Re\beta<1$.  One may then
define~\cite[Equation~1.4]{Ashyralyev2009}, \cite[Equation~(3.1)]{CarracedoAlix2001}
\[
D_0^\beta u=D_0^{\beta-1}D_0u=\frac{\sin\pi\beta}{\pi}\int_0^\infty
    \lambda^{\beta-1}(\lambda+D_0)^{-1}D_0u\,d\lambda
\quad\text{for $u\in X_1$,}
\]
and $D_0^\beta u=\I^{1-\beta}u'=(\I^{1-\beta}u)'$, that is, $D_0^\beta$ coincides with
both the Caputo and Riemann--Liouville fractional derivatives of order~$\beta$.

For $0\le\alpha\le1$, let $X_\alpha$ denote the complex interpolation
space~\cite[Chap.~4]{BerghLofstrom1976} arising from~$X_0$~and $X_1$.  The
next result is known~\cite[Theorem~11.6.1]{CarracedoAlix2001}, but
we sketch the proof in order to verify that the constant does not depend on~$\alpha$.

\begin{theorem}\label{thm:D alpha}
Let $0<\alpha<1$. For the spaces~\eqref{eq: X0 X1} and any complex Hilbert
space~$Z$, if the linear operator~$B$ satisfies $\|Bu\|_Z\le M_j\|u\|_{X_j}$ for $u\in X_j$ and $j\in\{0,1\}$, then
\[
\|Bu\|_Z\le e^{1+(\pi/4)^2}M_0^{1-\alpha}M_1^\alpha\|D_0^\alpha u\|_{X_0}
    \quad\text{for $D_0^\alpha u\in X_0$.}
\]
\end{theorem}
{\em Proof.}
If we set $Z_0=Z_1=Z$, then $Z_\alpha=Z$ with equal norms
\cite[Theorem~4.2.1]{BerghLofstrom1976}, so 
$\|Bu\|_Z\le M_0^{1-\alpha}M_1^\alpha\|u\|_{X_\alpha},$ and our task is to estimate the interpolation norm~$\|u\|_{X_\alpha}$ in terms of the fractional derivative~$D_0^\alpha u$.  Define the closed strip~$S=\{\,z\in\mathbb{C}:0\le\Re z\le1\,\}$ in the complex plane, and let $\mathcal{F}$ denote the space of functions $f:S\to X_0+X_1$ that are bounded and continuous on~$S$, analytic in the interior of~$S$,
and such that $f(iy)$ and $f(1+iy)$ tend to zero as~$|y|\to\infty$.  It can be
shown~\cite[Lemma~4.1.1]{BerghLofstrom1976} that $\mathcal{F}$ is a Banach
space with respect to the norm
\[
\|f\|_{\mathcal{F}}=\max\Bigl(\sup_{y\in\mathbb{R}}\|f(iy)\|_{X_0},
    \sup_{y\in\mathbb{R}}\|f(1+iy)\|_{X_1}\Bigr).
\]
The space~$X_\alpha$ then consists of those $u\in X_0+X_1$ such that
$u=f(\alpha)$ for some~$f\in\mathcal{F}$, with the interpolation norm defined by
\[
\|u\|_{X_\alpha}=\inf\{\,\|f\|_{\mathcal{F}}:
    \text{$u=f(\alpha)$ and $f\in\mathcal{F}$}\,\}.
\]
If we define $f(z)=e^{(z-\alpha)^2}D_0^{\alpha-z}u$, then $f(\alpha)=u$ with
$\|f(iy)\|_{X_0}=e^{\alpha^2-y^2}\|D_0^{-iy}D_0^\alpha u\|_{X_0}$ and 
\[\|f(1+iy)\|_{X_1}=e^{(1-\alpha)^2-y^2}\|D_0(D_0^{\alpha-1-iy}u)\|_{X_0}=e^{(1-\alpha)^2-y^2}\|D_0^{-iy}D_0^\alpha u\|_{X_0}.
\]
It is known~\cite[Example~2]{PrussSohr1990} that 
$\|D_0^{iy}u\|_{X_0}\le e^{\pi|y|/2}\|u\|_{X_0}$ for $u\in X_0$
because $D_0$ is m-accretive. Thus 
\[
\|u\|_{X_\alpha}\le\|f\|_{\mathcal{F}}
    \le\max\Bigl(\sup_{y\in\mathbb{R}}e^{\alpha^2-y^2+\pi|y|/2},
    \sup_{y\in\mathbb{R}}e^{(1-\alpha)^2-y^2+\pi|y|/2}\Bigr)
    \|D_0^\alpha u\|_{X_0},
\]
and since  $-y^2+\tfrac12\pi|y|
=(\tfrac14\pi)^2-(|y|-\tfrac14\pi)^2\le(\tfrac14\pi)^2,$ the proof is completed.
$\quad \Box$

The interpolation estimate used in the proof of Theorem~\ref{thm:Y} now follows
by repeating the preceding arguments with $D_0$~and $X_1$ replaced by their
time-reversed equivalents
\[
D_T u(t)=-u'(t)\quad\text{and}\quad
\XT{1}=\{\,u\in X_0:\text{$u'\in X_0$ and $u(T)=0$}\,\}.
\]
In fact, $\Re(u,D_T u)_{X_0}=\tfrac12\|u(0)\|_Y^2\ge0$ for all $u\in\XT{1}$,
and
\[
(\lambda+D_T)^{-1}u(t)=\int_t^Te^{-\lambda(s-t)}u(s)\,ds
    \quad\text{for $0\le t\le T$ and $\lambda\in\mathbb{C}$,}
\]
so $D_T$ is m-accretive. We find that if $0<\Re\beta<1$ then $D_T^{-\beta}u=\J{T}{\beta}u$, with 
\[
D_T^\beta u=D_TD_T^{\beta-1}u=(\J{T}{1-\beta}u)'
\quad\text{for $u\in\XT{1}$.}
\]

\begin{corollary}\label{cor:interpolation}
In the statement of Theorem~\ref{thm:D alpha}, we may replace $D_0$~and $X_1$
with $D_T$~and $\XT{1}$, respectively.
\end{corollary}


\end{document}